# Using Wavelets Based on B-splines for Calculation of the Hankel Transform

E.B. POSTNIKOV
Theoretical Physics Department,
Kursk State University
Radishcheva st, 33, Kursk, 305000
RUSSIA
postnicov@mail.ru

*Abstract:* The purpose of this paper is to represent the integral Hankel transform as a series. If one uses B-spline wavelet this series is a linear combination of the hypergeometrical functions. Numerical evaluation of the test function with known analytical Hankel transform of the null kind illustrates the proposed results.

## 1 Introduction

The Hankel transform [1] is a very useful tool of mathematical physics. Its application area is the axial symmetrical problems, especially seismology, electromagnetic methods in geophysical prospecting [2, 3], optical beam propagation [4] and other. There are two types of the Hankel transform. The first one is defined on the semi-infinite interval. In this case the direct and inverse transforms of the n-th kind are represented as an symmetric pair

$$F_n(p) = \int_0^\infty f(r) J_n(pr) r \, dr,$$
$$f(p) = \int_0^\infty F_n(p) J_n(pr) p \, dp. \quad (1)$$

In the case of the finite Hankel transform only a direct transform has an integral form. Without loss of generality its expression is

$$F_n(p) = \int_0^1 f(r) J_n(pr) r \, dr. \quad (2)$$

Unfortunately, comparatively little functions have a closed form solution of the Hankel integrals. But, there are several techniques for numerical evaluating of the Hankel transform. They are based on digital filtering, see a review in [2]; Discrere Hankel Transform (this is similar to the Fast Fourier Transform, see [5] and reference therein, and methods based on the Filon Quadratures. Filon's approach suggests that it is possible to fit only $f(r)$ (or $F_n(p)$) by a quadratic function instead of the entire integrand.

In this paper the generalization of this idea is considered. First note that one can use a more flexible basis as the intermediate agent in the calculation. The role of this intermediate agent can be played by the wavelet basis. This paper is organized as follows. In Section 2 we represent the integral Hankel transform as an exact series using B-spline wavelets of the arbitrary order. Section 3 consists of the numerical example of evaluation of the test function.

## 2 Outline of algorithm

Suppose there exist an exact analytical representation of the local weighted average of the function, and the weight is a power functions. Then there exist a representation of the integral Hankel transform as a series. The wavelet basis can play this role.

Suppose there exist a wavelet expansion of the function $f(x)$ [6]:

$$f(x) = \sum_{k=0}^\infty c_{ok} \varphi_k(r) + \sum_{j=0}^\infty \sum_{k=0}^\infty d_{jk} \psi_k(r),$$

where $\varphi_k(r) = \varphi(r-k)$, $\psi_{jk}(r) = 2^{j/2} \psi(2^j r - r)$ are a scaling function and a wavelet respectively.

Let us consider the semi-orthogonal wavelet basis based on the B-spline $N_n(x)$ [6]:

$$\psi(r) = \psi_m(r) = \sum_{n=0}^{3m-2} q_n N_n(2r-n). \quad (3)$$

This wavelet has a compact support $\operatorname{supp} \psi_m = [0, 2m-1]$, where $m$ is an order of the spline. The coefficients in (3) are

$$q_n = \frac{(-1)^n}{2^{m-1}} \sum_{l=0}^m \binom{m}{l} N_{2m}(n+1-l).$$

Denote the restriction of the B-spline to the unit interval by $N_m^\alpha$. Then, there exist an expansion of $N_m^\alpha$ into a set of Bernstein polynomials:

$$N_m^\alpha(y) = \sum_{l=0}^{m-1} a_l^{m-1}(\alpha) \binom{m-1}{\alpha} (1-[y-\alpha+1])^{m-l-1}(y-\alpha+1)^l$$

where $a_l^{m-1}(\alpha)$ are known Bernstein coefficients.

Take B-spline as a series

$$N_m(y) = \sum_{\alpha=1}^{m} N_m^\alpha(y) = \sum_{\alpha=1}^{m} N_m(y)\chi([\alpha-1,\alpha),y),$$

where $\chi([\alpha-1,\alpha),y) = \begin{cases} 1, \text{if } y \in [\alpha-1,\alpha) \\ 0, f\ y \notin [\alpha-1,\alpha) \end{cases}$.

Therefore, using Bernstein polynomials and a binomial expansion we get

$$\psi_m^{jk}(r) = \sum_{\alpha=1}^{m}\sum_{n=0}^{3m-2}\sum_{l=0}^{m-1}\sum_{\beta=0}^{m-l-1}\sum_{\gamma=0}^{\beta+l}(-1)^\beta \binom{m-1}{\alpha}\binom{m-l-1}{\beta}\binom{\beta+l}{\gamma} q_n a_l^{m-1}(\alpha)(1-2k-n-\alpha)^{\beta+l-\gamma} r^\gamma \times$$
$$\chi([\alpha+2k-n-1,\alpha+2k-n),r)$$

For any $r^\gamma$ there exist a definite integral [7]

$$\int_0^\zeta r^\gamma J_\nu(pr)rdr == \frac{p^\nu \zeta^{\gamma+2+\gamma}}{2^\nu(\gamma+2+\nu)\Gamma(\nu+1)} {}_1F_2\left(\begin{matrix}\frac{1}{2}(\gamma+2+\nu) \\ \frac{1}{2}(\gamma+4+\nu),\nu+1\end{matrix}\middle| -\frac{p^2\zeta^2}{4}\right)$$

Here ${}_1F_2$ is a hyperheometrical function.

Finally, the result of the Hankel transform of B-spline wavelet is the following analytical expression:

$$\psi_m^{jk}(p) = \sum_{\alpha=1}^{m}\sum_{n=0}^{3m-2}\sum_{l=0}^{m-1}\sum_{\beta=0}^{m-l-1}\sum_{\lambda=0}^{\beta+l}(-1)^\beta \binom{m-1}{\alpha}\binom{m-l-1}{\beta}\binom{\beta+l}{\gamma} q_n a_l^{m-1}(\alpha)(1-2k-n-\alpha)^{\beta+l-\gamma} \times$$
$$\frac{p^\nu}{2^{(j+2)\nu+2(j+1)}(\gamma+2+\nu)\Gamma(\nu+1)} \times$$
$$\left[(\alpha+2k-n)^{\gamma+2+\nu} {}_1F_2\left(\begin{matrix}\frac{1}{2}(\gamma+2+\nu) \\ \frac{1}{2}(\gamma+4+\nu),\nu+1\end{matrix}\middle| -\frac{p^2(\alpha+2k-n)^2}{2^{2(j+2)}}\right) - \right.$$
$$\left.(\alpha+2k-n-1)^{\gamma+2+\nu} {}_1F_2\left(\begin{matrix}\frac{1}{2}(\gamma+2+\nu) \\ \frac{1}{2}(\gamma+4+\nu),\nu+1\end{matrix}\middle| -\frac{p^2(\alpha+2k-n-1)^2}{2^{2(j+2)}}\right)\right] \quad (4)$$

Note that $c_{0k}$, $d_{jk}$ can be found by the convolution $f(r)$ with the dual basis functions. But these functions are also linear combinations of the power functions. This concludes the proof.

## 3 Numerical example
Consider a simple example [8]. The B-spline wavelet of the first kind is the Haar function:

$$\psi_1(r) = \begin{cases} 1, \ r \in (0,\tfrac{1}{2}) \\ -1, r \in (\tfrac{1}{2},1) \\ 0, \ r \notin (0,1) \end{cases}$$

The corresponding scaling function is

$$\varphi_1(r) = \begin{cases} 1, \ t \in (0,1) \\ 0, \ t \notin (0,1) \end{cases}$$

In this case (5) take a form

$$\psi_1^{jk}(p) = \frac{2^{-j}}{p}\left[2\left(k+\frac{1}{2}\right)J_1\left(2^{-j}p(k+\tfrac{1}{2})\right) - (k+1)J_1\left(2^{-j}p(k+1)\right) - kJ_1\left(2^{-j}pk\right)\right],$$

and the Hankel transform of the scaling function as

$$\varphi_1^{jk}(p) = \frac{1}{p}\left[J_1(p(k+1)) - J_1(p)\right].$$

The average and detail coefficients are

$$c_{0k} = \int_k^{k+1} f(r)dr,$$

$$d_{jk} = \int_{2^{-j}k}^{2^{-j}(k+1/2)} f(r)dr - \int_{2^{-j}(k+1/2)}^{2^{-j}(k+1)} f(r)dr.$$

We can compare he accuracy of the proposed algorithm. There is the known exact analytical Hankel transform:

$$\int_0^\infty e^{-\left(\frac{r}{a}\right)^2} J_0(pr) r\, dr = \frac{p}{2a} e^{-\left(\frac{p}{a}\right)^2}.$$

Let us truncate this transforming Gaussian function in the segment $[0, R]$. This function is represented in the Fig. 1 (the replacement $r$ to $x = r/h$ is used).

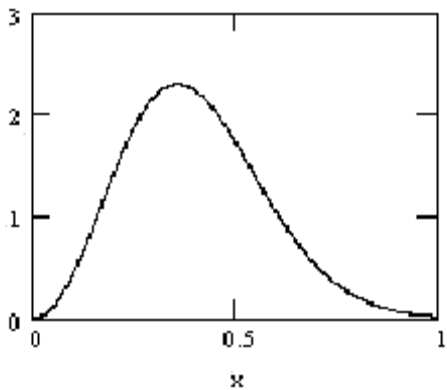

Fig .1

The exact transform (solid line) and approximate transform at the level of expansion $J = 3$ (dotted line) are plotted in Figure 2.

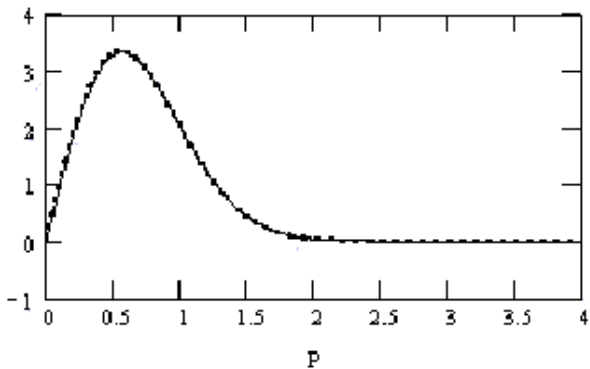

Fig. 2

The absolute error of this approximation is shown in Figure 3.

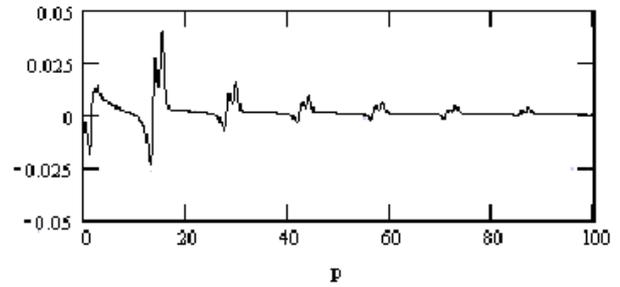

Fig. 3

One can see that the proposed method can be used for computational purposes.